\newcommand{\Z}{\mathbb{Z}}
\newcommand{\R}{\mathbb{R}}
\newcommand{\Sp}{\mathbb{S}}
\newcommand{\ind}{{\rm ind}}  
\mathchardef\varepsilon="010F
\mathchardef\epsilon="0122
\mathchardef\vartheta="0112
\mathchardef\theta="0123
\mathchardef\varrho="011A
\mathchardef\rho="0125
\mathchardef\varphi="011E     
\mathchardef\phi="0127
\renewcommand \emptyset \varnothing
\documentclass[12pt]{article}            
\usepackage[all]{xy}
\usepackage[latin1]{inputenc}        
\usepackage[dvips]{graphics}
\usepackage[dvips]{graphicx}
\usepackage{amsfonts}
\usepackage{amssymb}
\usepackage{amsmath}
\usepackage{amstext}
\usepackage{amsbsy}
\usepackage{amsopn}
\usepackage{amsthm}
\usepackage{amscd}
\usepackage{amsxtra}
\usepackage{mathrsfs}
\usepackage{upref}
\author{\sc Gianmarco Capitanio \\ 
{\it \small Universit{\' e} D. Diderot -- Paris VII, 
Equipe de G{\' e}om{\' e}trie et Dynamique} \\
{\it \small Case 7012 -- 2, place Jussieu , 75251 Paris Cedex 05} \\ 
{\small {\it e-mail:} {\rm Gianmarco.Capitanio@math.jussieu.fr}}
}
\date{}
\title{Generic singularities of minimax solutions to Hamilton--Jacobi
equations\footnote{To appear in {\it Journal of Geometry and
Physics.}} 
}

\begin{document}        

\numberwithin{equation}{section}                
\theoremstyle{plain}
\newtheorem{theorem*}{\bf Theorem}
\newtheorem{theorem}{\bf Theorem}
\newtheorem{lemma}{\bf Lemma}
\newtheorem{step}{\bf Step}
\newtheorem{proposition}{\bf Proposition}
\newtheorem*{proposition*}{\bf Proposition}
\newtheorem{corollary}{\bf Corollary}
\newtheorem*{corollary*}{\bf Corollary}
\theoremstyle{definition}
\newtheorem*{definition*}{\bf Definition}
\newtheorem*{definitions*}{\bf Definitions}
\newtheorem*{conjecture}{\bf Conjecture}
\newtheorem{example}{\bf Example}
\newtheorem*{example*}{\bf Example}
\theoremstyle{remark}
\newtheorem*{remark*}{\bf Remark}
\newtheorem{remark}{\bf Remark}
\newtheorem*{acknowledgements}{\bf Acknowledgements}

\maketitle

\begin{abstract}
Minimax solutions are weak solutions to Cauchy problems involving 
Hamilton--Jacobi equations, constructed from generating families 
quadratic at infinity of their geometric solutions.
We give a complete description of minimax solutions and we classify
their generic singularities of codimension not greater than $2$. 
\end{abstract}



\section{Introduction}\label{sct:1}

Hamilton--Jacobi equations play an important role in many fields of 
mathematics and physics, as for instance calculus of variations, optimal 
control theory, differential games, continuum mechanics and optics. 

Let us consider a Cauchy Problem involving a Hamilton--Jacobi equation.
For small enough time $t$ the solution $u$ is classically determined 
using the characteristic method. 
Although $u$ is initially smooth, there exists in general a critical time 
beyond which characteristics cross.
After this time, the solution $u$ is multivalued and singularities appear. 
Therefore, the problem how to extract ``true solutions''  
from multivalued solutions naturally arises. 

In 1991 Marc Chaperon proposed in \cite{Chaperon} a geometric method 
to construct weak solutions to Hamilton--Jacobi equations, called
minimax solutions (see also \cite{Viterbo2}).
Their definition is based on generating families quadratic at infinity of 
Lagrangian submanifolds. 
The minimax solutions have the same analytic properties as the 
viscosity solutions, namely existence and uniqueness theorems hold. 
However, minimax solutions are in general different from viscosity
solutions. 

This paper is organized as follows. 
The first part is a survey on minimax solutions. 
In the second part we study the small codimension 
generic singularities of minimax solutions.
Namely, we prove that the singularities of codimension not greater
than $2$ of minimax solutions and viscosity solutions are the same. 
This result may be useful to find examples in which minimax 
solutions have physical meaning. 

The survey part of the paper is mainly based on a talk I have given 
for the ``Trimester on Dynamical and Control Systems'' at SISSA-ICTP 
in Trieste. 
I thank Andrei Agrachev for its kind invitation. 

\section{Minimax of functions quadratic at infinity}\label{sct:2}

Let $E$ be the space $\R^n$ and $f:E\rightarrow \R$ a smooth 
function having a finite number if critical points.  
We suppose $f$ {\it quadratic at infinity}, that is
$f$ is a non-degenerate quadratic form outside a compact set.
For $\lambda\in\R$, we set 
$$E^\lambda:= \{\xi\in E : f(\xi)\leq \lambda\}\ . $$
Let us fix $a\in\R$, big enough for that $(-a,a)$ contains all the
critical values of $f$. 
Then we set $E^{\pm\infty}:=E^{\pm a}$. 

Let us denote by $\rho_\lambda$ the restriction mapping 
$E\rightarrow E^\lambda$, inducing the homomorphism 
$\rho^*_\lambda: \tilde H_*(E,E^{-\infty}) \rightarrow \tilde
H_*(E^\lambda,E^{-\infty})$ 
between the relative reduced $\Z$-homology groups. 
The quotient space $E/E^{-\infty}$ is homeomorphic to $\Sp^q$, where 
$q$ denotes the index of the quadratic part of $f$ at infinity. 
Hence, the only non-trivial reduced homology group is $ \tilde
H_q(E,E^{-\infty})\simeq \Z$.  
Fix a generator $\gamma$ of this group. 

\begin{definition*}
The {\it minimax} of $f$ is the real number  
$$\min\max (f) := \inf\{\lambda\in\R : \rho^*_\lambda\gamma\not=0\}\ .$$
\end{definition*} 

The infimum defining the minimax exists and it is
finite, indeed $\rho^*_\lambda\gamma=\gamma\not=0$ and 
$\rho^*_{-\lambda}\gamma=0$ whenever $\lambda>a$; it does not depend on 
the choice of the generator $\gamma$. 
Moreover, the minimax of a function is a critical value, since by 
definition the topology of the sublevel sets $E^\lambda$ changes when
$\lambda$ crosses the minimax value. 

Since the minimax level is stable under small
deformation of the function (see \cite{io2003}), we may assume without
loss of generality that $f$ is generic. 
Genericity conditions mean that all the critical points of $f$ are of
Morse type and all its critical values are different. 
Therefore, every critical point has an {\it index}, which is
the number of negative signs in the Morse normal form of the function
at the critical point. 

Morse Theory (see \cite{milnor}) describes how the topology of the
sublevels sets $E^\lambda$ changes when $\lambda$ changes. 
Namely, {\it $E^\lambda$ is diffeomorphic to $E^\mu$, provided that 
$[\lambda,\mu]$ does not contain critical values}. 
On the other hand, {\it if $[\lambda,\mu]$ contains only one critical 
value, then $E^\mu$ retracts on the space obtained from $E^\lambda$ 
attaching to its boundary a cell of dimension equal to the index of 
the critical point realizing the critical value}. 
Hence, every generic function on $E$ defines a cell decomposition
of the space and there is a $1:1$ correspondence between its cells and
the critical points of the function.  

Let us recall that the {\it incidence coefficient} of two cells in this 
decomposition is the degree of the restriction of the attaching map 
from the higher dimension cell to the lower dimension cell 
(see \cite{DNF}). 
The {\it incidence coefficient} $[\xi:\eta]$ of two critical points $\xi$ and
$\eta$ of $f$ is the incidence coefficient of the corresponding cells. 
Note that  $[\xi:\eta]\not=0$ implies that $\ind(\xi)=\ind(\eta)+1$
and $f(\xi)>f(\eta)$. 
For example, the two critical points of $f(\xi)=\xi^3-\epsilon \xi$ have
non-zero incidence coefficient for $\epsilon>0$. 

The incidence coefficient of pairs of critical points leads to a natural
partition of the critical set $\Sigma$ of $f$ in the following way. 
Denote by $C_1$ a pair of critical points realizing the minimum of the
set 
$$\{f(\xi)-f(\eta) : \xi,\eta\in\Sigma, [\xi:\eta]\not=0 \} \ .$$
Then define by induction $C_{i+1}$ from $C_1,\dots,C_i$ as a pair of
critical points realizing the minimum of the set 
$$\{f(\xi)-f(\eta) : \xi,\eta\in\Sigma\smallsetminus (C_1\cup \dots\cup C_i) , 
[\xi:\eta]\not=0 \} \ .$$
In this way, we decompose the critical set into the
disjoint union of pairs $C_i$ and a set $F$ which does not contains
incident critical points, i. e. $[\xi:\eta]=0$ for every $\xi,\eta\in F$. 

\begin{definition*} 
Two critical points are {\it coupled} if they belong to the
same pair $C_i$ in the preceding decomposition;  
a critical point is {\it free} if it belongs to $F$.
\end{definition*}

In \cite{io2003} we proved the following result.

\begin{theorem}\label{thm:1} 
Every generic function quadratic at infinity has exactly one free
critical point, and its value is the minimax. 
\end{theorem}

\section{Geometric and multivalued solutions}\label{sct:3}

In this section we introduce the symplectic framework for the
geometric and multivalued solutions 
of Hamilton--Jacobi equations. 
Let $X$ be a closed manifold of dimension $n$, $\pi:T^*X\rightarrow X$
its cotangent bundle, endowed with the standard symplectic form
$d\lambda$, where $\lambda$ is the Liouville's $1$-form. 
In local coordinates $T^*X=\{x,y\}$, we have $\pi(x,y)=x$ and
$\lambda=y\ dx$. 
A dimension $n$ submanifold of $T^*X$ is said to be {\it Lagrangian}
whenever $d\lambda$ vanishes on it. 
Two Lagrangian submanifolds are {\it Hamiltonian isotopic}  
if there exists a flow generated by a Hamiltonian field, transforming
one into the other.   
Hamiltonian isotopies transform Lagrangian submanifolds into
Lagrangian submanifolds. 
A {\it generating family} of a Lagrangian submanifold $L\subset T^*X$ is a 
smooth function $S: X\times \R^k\rightarrow\R$
such that
$$L = \big\{ \left(x,\partial_x S (x;\xi)\right) : 
       \partial_\xi S (x;\xi)=0 \big\} \ ,$$
where $S$ verifies also the rank condition 
${\rm rk} \left( \partial^2_{\xi\xi} S,\partial^2_{\xi x} 
S \right)|_{\partial_\xi S=0} \ = \max$.

Given a generating family $S$, the following operations give rise to 
new generating families $T$ of the same Lagrangian submanifold: 
\begin{enumerate}
\item {\it Addition of a constant:} 
$T(x;\xi)=S(x;\xi) + C$ for $C\in\R$; 
\item {\it Stabilization:} $T(x;\xi, \eta)=S(x;\xi)+Q(\eta)$, where 
$Q$ is a non-degenerate quadratic form;  
\item {\it Diffeomorphism:} 
$T(x;\eta)=S(x;\xi(x,\eta))$, where $(x;\eta) \mapsto
(x,\xi(x,\eta))$ is a fibered diffeomorphism with respect to the
coordinate $x$. 
\end{enumerate} 

Two generating families are {\it equivalent} if one can obtain one
from the other by a finite sequence of the preceding operations. 
A generating family $S$ is {\it quadratic at infinity} (gfqi) if there exists a 
non-degenerate quadratic form $Q$ such that 
$S(x;\xi)=Q(\xi)$, whenever $|\xi|$ is big enough (uniformly in $x$).

\begin{theorem}[Sikorav--Viterbo, \cite{Sikorav}---\cite{Viterbo}] 
Every Lagrangian submanifold of $T^*X$, Hamiltonian isotopic to the zero section 
$\{(x;0) :  x\in X\}$, admits a unique gfqi modulo the preceding
equivalence relation.
\end{theorem}

\begin{remark*}
The Theorem still holds in the case of non-compact manifolds,
provided that the projection of the Lagrangian submanifold into the
base is $1:1$ outside a compact set. 
\end{remark*}

A Lagrangian submanifold $L\subset T^*X$ is {\it exact} if the
Liouville's $1$-form is exact on it. 
In this case, we can lift $L$ to a Legendrian submanifold in the
contact space $J^1X$ of $1$-jets over $X$. 
This lifting projects into a wave front in the space
$J^0X\simeq X\times\R=\{x,z\}$ of $0$-jets over $X$, defined
up to shifts in the $z$-direction.  
If $S$ is a gfqi of $L$, then the corresponding wave front is the
graph of $S$ 
$$\left\{ (x,S(x;\xi)) : \partial_\xi S(x;\xi)=0\right\} \ .$$
Graphs of equivalent gfqi are equal up to a translation
along the $z$-axis. 
The points for which the natural projection ``forgetting $z$'' of the
graph into $X$ is not a fibration form a stratified hypersurface in
$X$, provided that $L$ is generic.

We may now define generalized solutions to Hamilton--Jacobi equations. 
Let $Q$ be a closed manifold. 
We consider the following Cauchy Problem on $Q$ involving a Hamilton--Jacobi
equation:
$$(CP)\begin{cases}
\partial_t u(t,q) + H(t,q,D_q u(t,q))=0 , &\forall \ t\in\R^+ , q\in Q\\
u(0,q)=u_0(q) , &\forall\ q \in Q \ .
\end{cases}$$
The Hamiltonian $H:\bar \R^+\times T^*Q\rightarrow \R$ is supposed of 
class $\mathscr C^2$ on $\R^+\times T^*Q$ and continuous at the
boundary; the initial condition $u_0:Q\rightarrow\R$ is assumed to
be of class $\mathscr C^1$.

Let $M:=Q\times\R$ be the space-time manifold and 
$T^*M =\{t,q;\tau,p\}$ its cotangent bundle (endowed with the 
symplectic form $dp\wedge dq+d\tau\wedge dt$). 
We consider the flow 
$\Phi:\bar \R^+\times T^*M\rightarrow  T^*M$,   
generated by the Hamiltonian $\mathscr H:T^*M\rightarrow\R$ defined
by $\mathscr H:=\tau+H$, and the submanifold
$$\sigma:=\left\{ \big((0,q;-H(0,q, du_0(q)),du_0(q) \big) : q\in
Q\right\}\subset  \mathscr H^{-1}(0)\subset T^*M\ .   $$

\begin{definition*}
The {\it geometric solution} of $(CP)$ is the submanifold 
$$L:=\bigcup_{t>0} \Phi^t(\sigma) \subset T^*M \ .$$
\end{definition*}

It turns out that the geometric solution to $(CP)$ is an exact
Lagrangian submanifold, contained into the hypersurface $\mathscr
H^{-1}(0)$ and Hamiltonian isotopic (for finite times) to the zero
section of $T^*M$ (see \cite{io2003}).
Hence, the Sikorav--Viterbo Theorem guarantees that it has a unique gfqi
modulo the equivalence relation.  
Up to a suitable constant, we may assume that the graph of
the gfqi restrained at $t=0$ is equal to $\sigma$. 
We denote by $S$ such a gfqi. 

Let us consider now the projection ${\rm pr}:T^*M\rightarrow T^*Q$,
defined locally by 
${\rm pr}(t,q;\tau,p):=(q,p)$.
We call {\it geometric solution at time $t$} of $(CP)$ the exact
Lagrangian submanifold 
$L_t:= {\rm pr}\circ\Phi^t(\sigma)\subset T^*Q$,
which is an isochrone section of the global geometric solution.  
It is easy to check that $L_t$ is Hamiltonian isotopic to the zero
section, so, by the Uniqueness Theorem, its gfqi is
$S_t(q;\xi):=S(t,q;\xi)$.  

\begin{definition*}
The graph of $S$ (resp., $S_t$) will be called the {\it multivalued
  solution} (resp., {\it at time $t$}) to $(CP)$.
\end{definition*} 

\begin{remark*}
It is possible to construct global generating families of geometric 
solutions as follows. 
The  {\it action functional} $\int p dq - H dt$ is a global generating
family, whose parameters belong to an infinite dimension space.   
By a fix point method, proposed by Amann, Conley and Zehnder, one
can obtain a true generating family (with finite dimensional
parameters), 
see \cite{Cardin}. 
\end{remark*}

\section{Minimax and viscosity solutions} 

In this section we introduce minimax solutions to Hamilton--Jacobi
equations, and we discuss their properties and relations with
viscosity solutions. 
Let $S$ be the gfqi of the geometric solution to $(CP)$. 
For every fixed point $(t,q)$ in the space time, 
the function $S_{t,q}$ is quadratic at infinity, so we
may consider its minimax critical values, defined in section \ref{sct:2}. 

\begin{definition*}
The {\it minimax solution} to $(CP)$ is the function defined by 
$$u(t,q):= \min\max (S_{t,q}) \ . $$
\end{definition*}

\begin{remark*}
For a Cauchy Problem $(CP)$ on a non-compact manifold, as for instance
for $Q=\R^n$, we define the minimax solution as follows. 
Let $(H_n)_n$ be a sequence of smooth Hamiltonians with compact support
in the variables $q$. 
Then, at any Cauchy Problem defined by $H_n$ and the
given initial condition, we may associate its minimax solution $U_n$,
defined as above. 
By definition, the minimax solution of $(CP)$ is the limit of the
solutions $U_n$ for $n\rightarrow \infty$.
\end{remark*}

\begin{theorem}[Chaperon, \cite{Chaperon}] 
The minimax solutions  are  weak solutions\footnote{A weak solution to
  $(CP)$ is a continuous and almost everywhere derivable function,
  which solves in these points the Hamilton--Jacobi equation;
  moreover, its restriction at $t=0$ is equal to the initial
  condition.}  
to $(CP)$, Lipschitz on finite times, which not depend on the choice of the gfqi.
\end{theorem}

We end this section discussing the relations between minimax and
viscosity solutions.  
Viscosity solutions to Hamilton--Jacobi equations have been introduced by 
Crandall and Lions in \cite{crandall}, following the generalization
proposed by Kruzhkov of the Hopf's formulas. 

\begin{theorem}[Joukovskaia, \cite{Joukovskaia}]
The minimax and the viscosity solutions of $(CP)$ are equal, provided that the
Hamiltonian defining the equation is convex or concave in the
variables $p$.
\end{theorem}

Izumiya and Kossioris constructed in \cite{izumiya} the viscosity
solution beyond its first critical time. 
Actually, they proved that the Lagrangian graph of the viscosity
solution is not always contained into the geometric solution of the equation. 
Hence, minimax and viscosity solutions are in general different. 

\section{Characterization of minimax solutions} 

In order to describe a useful characterization of minimax solutions 
(first presented in \cite{io2003}), we start introducing some standard 
notations concerning unidimensional wave fronts. 

\begin{figure}[h]    
  \begin{center}
  \scalebox{.3}{\input{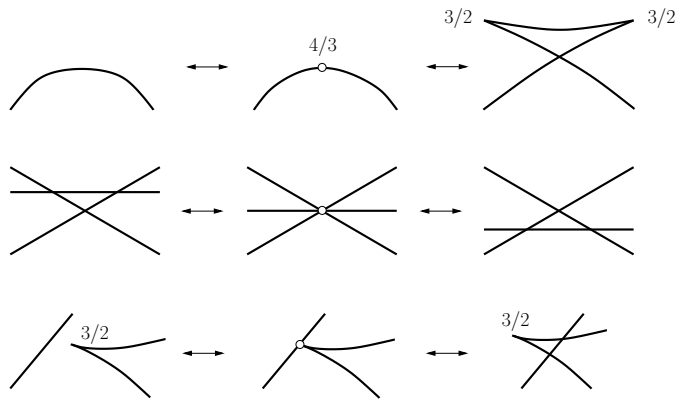}}
  \end{center}
  \caption{Generic singularities of flat fronts under Legendrian  isotopies.}  
\label{fig:4}
\end{figure}

Let $J^0\R\simeq\R^2=\{(q,z)\}$ be the space of $0$-jets over $\R$ and 
$\pi_0 : J^0\R\rightarrow\R$ the natural fibration $(q,z)\mapsto q$. 
A {\it wave front} in $J^0\R$ is the projection of an embedded 
Legendrian curve in the contact space $J^1\R\simeq\R^3=\{(q,z,p)\}$ 
under the projection $\pi_1 : (q,z,p) \mapsto (q,z)$. 
The only singularities of generic wave fronts are semicubic
cusps and transversal self-intersections.
Two embedded Legendrian curves in $J^1\R$ are {\it Legendrian
isotopic} if there exists a smooth path joining them in the space 
of the embedded Legendrian curves. 
In this case their projections are also called isotopic. 
A front is {\it long} if it is the graph of a smooth function outside
a compact set of $\R$; it is {\it flat} if its tangent lines are
nowhere vertical. 
A {\it section} of a flat front is a connected maximal subset which 
is the graph of a piecewise $\mathscr C^1$ function; 
a {\it branch} is a smooth section. 
The generic perestroikas occurring to flat wave fronts under 
isotopies among flat fronts are the following, illustrated in figure \ref{fig:4}: 
{\it cusp birth/death}, {\it triple intersection}, 
{\it cusp crossing}.

We come back now to the Cauchy Problem $(CP)$ in the case $Q=\R$. 
We assume that its geometric solution is a generic Lagrangian
curve. 
A multivalued solution $\Sigma$ of $(CP)$ at a generic time is a flat
wave front, isotopic to $\{(q,0):q\in\R\}$.  
Up to replace our Cauchy Problem with a sequence of problems
approximating it, we may assume also that $\Sigma$ is a long wave
front, and that it is the graph of a gfqi, denoted by $S$, of the geometric
solution.  
Hence, every branch of the multivalued solution is realized by 
a critical point of $S_{t,q}$ (depending on $q$ as a parameter).   

\begin{definition*}
The {\it index} of a branch is by definition the index of its
corresponding critical point minus the index of the quadratic part of
$S$. 
\end{definition*}
  
The index of a branch does not depend on the choice of the gfqi. 
The non-compact branch of $\Sigma$ has index $0$, as well 
as the minimax section.
Moreover, the index changes by $+1$ (resp., $-1$) passing through a
positive cusp\footnote{A cusp is {\it positive} (resp., {\it negative}) if,
  following the orientation of the front, we pass from a branch to the
  other according to the coorientation defined by the positive $z$
  direction.} 
(resp., negative cusp).  
A double point of $\Sigma$ is called {\it homogeneous} whenever it is
the intersection of two branches having the same index. 

The partition of the critical set of a generic function quadratic at
infinity we give in section \ref{sct:2} leads, by Theorem \ref{thm:1}, 
to a decomposition of the multivalued solution 
into the graph $\mu$ of the minimax solution
and the pairs $X_i$ of coupled sections, formed by the sections 
corresponding to coupled critical points of the gfqi:
$$\Sigma=\mu \cup \bigcup_i X_i \ .$$ 
Indeed, when $q$ runs on $\R$, every pair of coupled critical points of
$S_{t,q}$ moves on the front along two sections, forming a closed
curve.   

The curves $X_i$ are homeomorphic to a disc boundary, they have
exactly two cusps and no self-intersections.
They are almost everywhere smooth and only continuous at homogeneous double
points of the multivalued solution (see \cite{io2003}). 

\begin{figure}[h]    
  \begin{center}
  \scalebox{.3}{\input{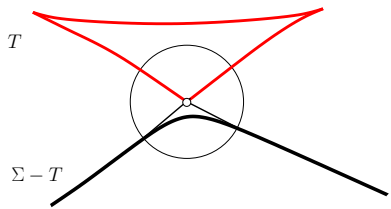}}
  \end{center}
  \caption{The front $\Sigma-T$.}  
\label{fig:5}
\end{figure}

Any homogeneous double point of an oriented multivalued solution 
defines in it a closed connected curve. 
Such a curve is called a {\it triangle} whenever it has exactly two cusps.  
The double point is the {\it vertex} of the triangle.  
Let $T$ be a triangle of $\Sigma$. 
Fix an arbitrary small ball centered at the vertex of $T$,
intersecting only the two branches containing the double point. 
We denote by $\Sigma-T$ a long flat front equal to
$\Sigma\smallsetminus T$ outside the ball and smooth inside it 
(see figure \ref{fig:5}).

\begin{definition*}
A triangle $T$ of $\Sigma$ is {\it vanishing} if there exists a 
Legendrian isotopy joining  $\Sigma$ and  $\Sigma-T$ among long flat
fronts without homogeneous triple intersections. 
\end{definition*}
 
The following Theorem (see \cite{io2003}) provides an effective
method to simplify recursively a given multivalued solution.  
After a finite number of steps, we get a front, which is actually the
graph of the minimax solution to the initial problem.

\begin{theorem}\label{thm:6}
Let $\Sigma$ be a multivalued solution and 
$\mu\cup X_1 \cup\dots\cup X_N$ its decomposition into minimax section
and coupled sections. 
Then $\Sigma$ is smooth or has a vanishing triangle $T$ among the
curves $X_i$.  
In this case, the decomposition of $\Sigma-T$ is induced by that of $\Sigma$. 
In particular, the minimax of $\Sigma$ and $\Sigma-T$ are equal outside 
the ball centered at the vertex of the vanishing triangle.
\end{theorem}

\begin{example*}
Consider the multivalued solution $\Sigma$ depicted in figure \ref{fig:6}. 
The homogeneous double points $B$, $C$, $D$, $F$, $G$ and $I$ are not verteces of
triangles.  
The triangles $T_E$ and $T_L$ of verteces $E$ and $L$ are vanishing, 
while those of verteces $A$, $H$ and $M$ are not. 
Indeed, to cut off $T_A$ and $T_M$ we must pass through a
self-tangency, while to cut off $T_H$ we must pass through a homogeneous triple
intersection at $F$.  
The minimax of $\Sigma$ is equal, outside two arbitrary small balls
centered at $E$ and $L$, to the minimax of $\Sigma':=\Sigma-T_E-T_L$. 
Now, $B$ and $I$ define two vanishing triangles $T_B$ and $T_I$; 
$\Sigma'- T_B-T_I$ is smooth, so it is the minimax section of
$\Sigma$ outside two arbitrary small balls centered
at $B$ and $I$.
\begin{figure}[h]    
  \begin{center}
  \scalebox{.39}{\input{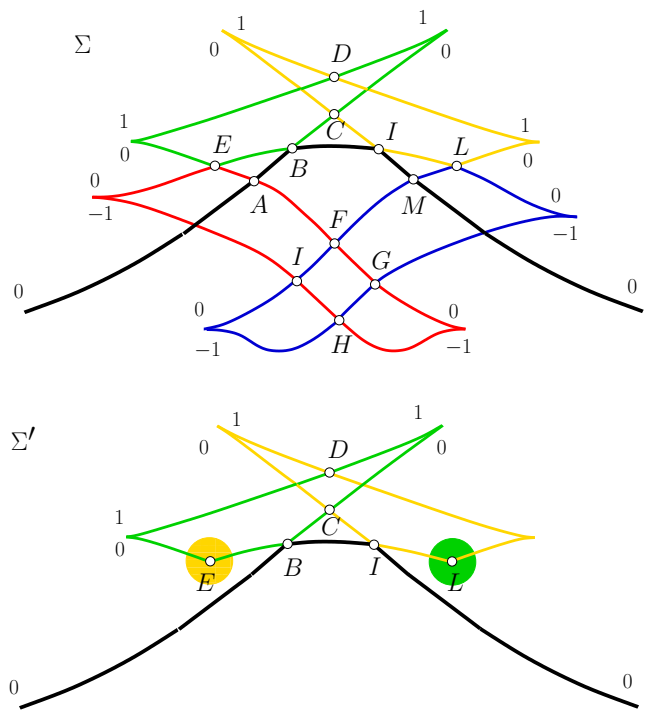}}
  \end{center}
  \caption{The multivalued solutions $\Sigma$ and $\Sigma'$.}  
\label{fig:6}
\end{figure}
\end{example*}

\section{Singularities of minimax solutions}\label{sct:8}

In this section we classify the generic singularities of small
codimension arising in minimax solutions to Hamilton--Jacobi
equations. 
This classification is made with respect to the Left--Right
equivalence, fibered on the time direction, defined below. 
Let us consider two functions $f,g:\R\times\R^n\rightarrow \R$. 
We denote by $t$ the coordinate in $\R$ and
by $q=(q_1,\dots,q_n)$ the coordinate in $\R^n$. 
 
\begin{definition*}{\rm 
Let $a,b$ be two points of the space-time $\R\times\R^n$.  
The germ of $f$ at $a$ and the germ of $g$ at $b$ are said to be 
{\it equivalent} if there exists two diffeomorphism germs $\phi$ and
$\psi$ such that $\psi\circ f\circ \phi^{-1}=g$ and $\phi$ is fibered
with respect to the time axis:
$\phi(t,q)=\big(T(t),Q(t,q)\big)$. 
}\end{definition*}

A {\it singularity} of a function germ for this equivalence relation 
is its equivalence class. 
A Cauchy Problem on $Q=\R^n$ is said to be {\it generic} if its geometric
solution is generic as Lagrangian submanifold; its minimax solution 
is also called {\it generic}. 

In \cite{Joukovskaia},  Joukovskaia proved that the singular set
of any generic minimax solution (formed by the points 
where the solution is not $\mathscr C^1$) is a closed stratified
hypersurface, diffeomorphic at any point to a semi-algebraic
hypersurface.  

\begin{theorem}\label{thm:8}
Let $u$ be a generic minimax solution and $(t,q)$ a point belonging to a
stratum of codimension $c=1$ or $2$ of its singular set.  
Then the germ of $u$ at $(t,q)$ is equivalent to one of the map germs
in the table below. 
\begin{center}
\begin{tabular}{|| c | c ||} 
\hline \hline
$c$ & Normal form\\ \hline\hline
$1$ & $|q_1|$ \\ \hline 
$2$ & $ \min \{|q_1|,t\} $ \\ \hline  
$2$ & $\min \{Y^4-tY^2+q_1Y : Y\in\R\} $ \\ \hline 
\end{tabular}
\bigskip\end{center}
\end{theorem}

\begin{figure}[h]
  \centering
   \scalebox{.45}{\input{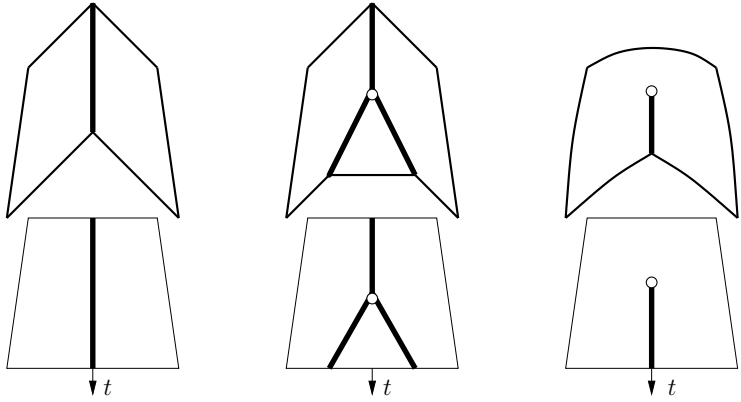}}
  \caption{Small codimension generic singularities of minimax solutions.}
 \label{fig:9}
\end{figure}

The images of these normal forms (in the case $n=2$) are depicted 
in figure \ref{fig:9}.
Note that these singularities are stable.
These map germs are also the normal forms of the generic singularities (of
codimension not greater than $2$) of viscosity solutions 
(whose classification has been done by Bogaevski in \cite{Bogaevski}).

Joukovskaia studied in \cite{Joukovskaia} the generic singularities of
minimax functions, defined as minimax levels of Lagrangian
submanifold's gfqi which are not necessarily geometric solutions of
Hamilton--Jacobi equations. 
To prove Theorem \ref{thm:8}, we will show that some singularities in
Joukovskaia's list do not arise in minimax solutions. 
In order to do this, let us recall her Theorem. 

\begin{theorem}[Joukovskaia, \cite{Joukovskaia}]
Any generic minimax function on $\R\times \R^n$ is equivalent, at
any point in a codimension $2$ stratum of
its singular set, to one of the map germs listed in Theorem
\ref{thm:8} (for $c=2$) or
to one of the following map germs:  
\begin{enumerate}
\item[{\rm (a)}] $(t;q_1,\dots,q_n)\mapsto \max\big\{t,-|q_1|\big\}$;
\item[{\rm (b)}] $(t;q_1,\dots,q_n)\mapsto \min\big\{|q_1|,\max\{-|q_1|,t\}\big\}$.
\end{enumerate}
\end{theorem}
 
\begin{figure}[htpb]
  \centering
   \scalebox{.42}{\input{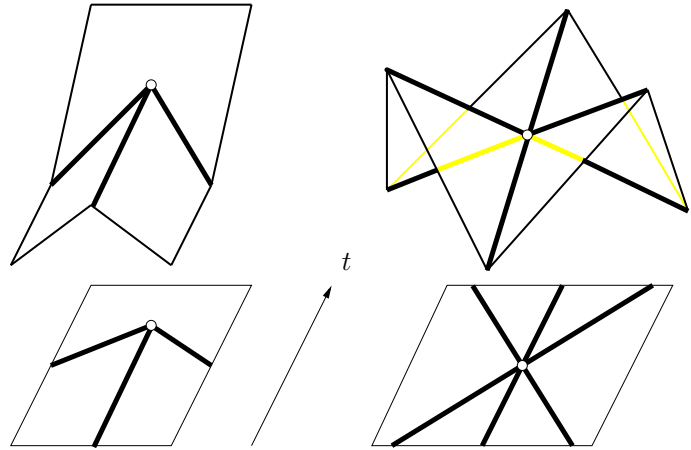}}
  \caption{Further generic singularities of minimax functions.}
 \label{fig:10}
\end{figure}

The images of the map germs (a) and (b) are shown in figure \ref{fig:10}.

\begin{proof}[Proof of Theorem \ref{thm:8}]  
We may assume $n=1$ without loss of generality.
Let us consider a multivalued solution having a singularity (a) or
(b).  
In both cases, there are three $0$-index branches crossing at the
singular point (settled at the origin). 

Fix an arbitrary small neighborhood $U$ of the origin in $\R$;
then for every small enough  time $t<0$, the isochrone 
multivalued solution at this time takes, over the corresponding 
section of $U$, the configuration illustrated in the left
part of figure \ref{fig:11}; moreover, all the other sections of the
front's decomposition are smooth over $U$.   
The minimax section $\mu$ has three jumps in $A$, $B$ and $C$, 
shrinking to the homogeneous triple point at the origin as $t$ goes to
$0$. 

Actually, the jumps $A$ and $C$ belong to the same section, denoted by
$\alpha$ in the figure.   
\begin{figure}[htpb]
  \centering
   \scalebox{.33}{\input{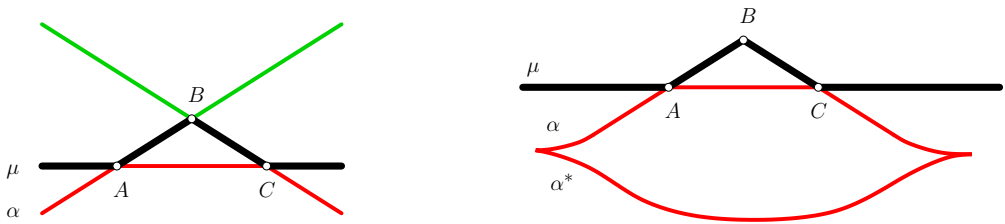}}
  \caption{Multivalued solution configurations near the singular point.}
 \label{fig:11}
\end{figure} 
By Theorem \ref{thm:6} we may recursively eliminate the 
vanishing triangles in the multivalued solution changing neither the
minimax section nor the other pairs of coupled critical points. 
Hence, after a finite number of such operations, the front is composed
by the minimax solutions, the section $\alpha$ and its
coupled section $\alpha^*$, see the right part of figure \ref{fig:11}. 

This front is not isotopic to the trivial front $\{(q,0):q\in\R\}$. 
Hence, multivalued solutions have neither singular points of type 
$(a)$ nor $(b)$.
Theorem \ref{thm:8} is now proven.
\end{proof}

\begin{remark*}
In other words, both  singularities (a) and (b) imply the existence of a
cycle in the graph associated to a decomposition of an isochrone
multivalued solution of the Cauchy Problem. 
This is impossible, since all the graphs associated to admissible
decompositions of isochrone multivalued solutions are trees, due to
Theorem \ref{thm:6}.  
\end{remark*}


\end{document}